\documentclass[11pt]{article}
\usepackage{enumerate}
\usepackage{amssymb,a4wide,latexsym,makeidx,epsfig,fleqn}
\usepackage{amsthm}
\usepackage{amsmath}
\usepackage{enumerate}
\newtheorem{theorem}{Theorem}[section]

\newtheorem{definition}{Definition}
\newtheorem{lemma}[theorem]{Lemma}

\begin{document}
\textwidth 150mm \textheight 225mm
\title{Complex unit gain bicyclic graphs with rank 2, 3 or 4
\thanks{ Supported by the National Natural Science Foundation of China (No. 11171273).}}
\author{{Yong Lu, Ligong Wang\footnote{Corresponding author} and Peng Xiao}\\
{\small Department of Applied Mathematics, School of Science, Northwestern
Polytechnical University,}\\ {\small  Xi'an, Shaanxi 710072,
People's Republic
of China.} \\{\small E-mails: luyong.gougou@163.com; lgwangmath@163.com; xiaopeng@sust.edu.cn}}

\date{}
\maketitle
\begin{center}
\begin{minipage}{120mm}
\vskip 0.3cm
\begin{center}
{\small {\bf Abstract}}
\end{center}
{\small A $\mathbb{T}$-gain graph is a triple $\Phi=(G,\mathbb{T},\varphi)$ consisting of a graph $G=(V,E)$, the circle group $\mathbb{T}=\{z\in C: |z|=1\}$ and a gain function $\varphi:\overrightarrow{E}\rightarrow \mathbb{T}$ such that $\varphi(e_{ij})=\varphi(e_{ji})^{-1}=\overline{\varphi(e_{ji})}$. The rank of $\mathbb{T}$-gain graph $\Phi$, denoted by $r(\Phi)$, is the rank of the adjacency matrix of $\Phi$. In 2015, Yu, Qu and Tu [ G. H. Yu, H. Qu, J. H. Tu, Inertia of complex unit gain graphs, Appl. Math. Comput.  265(2015) 619--629 ] obtained some properties of inertia of a $\mathbb{T}$-gain graph. They characterized the $\mathbb{T}$-gain unicyclic graphs with small positive or negative index. Motivated by above, in this paper, we characterize the complex unit gain bicyclic graphs with rank 2, 3 or 4.

 }
\end{minipage}
\end{center}

\section{Introduction }
All graphs considered in this article are simple graphs. Let $G=(V,E)$ be a simple graph with vertex set $V=V(G)$ and edge set $E=E(G)$. A \emph{gain graph} is a graph whose edges are labeled orientably by elements of a group $M$. That is, if an edge $e$ in one direction has label a group element $m$ in M, then in the other direction it has label $m^{-1}$ (the invertible element of $m$ in $M$). We call the group $M$ to be the \emph{gain group}. A gain graph is a generalization of a signed graph, where the gain group $M$ has only two elements 1 and $-1$, see Zaslavsky \cite{ZBI}.

 A $\mathbb{T}$-\emph{gain graph (or complex unit gain graph)} is a graph with the additional structure that each orientation of an edge is given a complex unit, called a \emph{gain}, which is the inverse of the complex unit assigned to the opposite orientation. For a simple graph $G=(V,E)$ of order $n$, let $\overrightarrow{E}$ be the set of oriented edges, it is obvious that this set contains two copies of each edge with opposite directions. We write $e_{ij}$ for the oriented edge from $v_{i}$ to $v_{j}$. The \emph{circle group}, which is  denoted by $\mathbb{T}=\{z\in C: |z|=1\}$, is a subgroup of the multiplicative group of all nonzero complex numbers $\mathbb{C}^{\times}$.  A $\mathbb{T}$-gain graph is a triple $\Phi=(G,\mathbb{T},\varphi)$ consisting of a graph $G=(V,E)$, the circle group $\mathbb{T}=\{z\in C: |z|=1\}$ and a gain function $\varphi:\overrightarrow{E}\rightarrow \mathbb{T}$ such that $\varphi(e_{ij})=\varphi(e_{ji})^{-1}=\overline{\varphi(e_{ji})}$, where $G$ is the \emph{underlying graph} of the $\mathbb{T}$-gain graph. We often write $\Phi=(G,\varphi)$ or $G^{\varphi}$ for a $\mathbb{T}$-gain graph. The \emph{adjacency matrix} associated to the $\mathbb{T}$-gain graph $\Phi$ is the $n\times n$ complex matrix $A(\Phi)=(a_{ij})$, where $a_{ij}=\varphi(e_{ij})$ if $v_{i}$ is adjacent to $v_{j}$, otherwise $a_{ij}=0$. It is obvious to see that $A(\Phi)$ is Hermitian and its eigenvalues are real. If the gain of every edge is 1 in $\Phi$, then the adjacency matrix $A(\Phi)$ is exactly the adjacency matrix $A(G)$ of the underlying graph $G$. It is obvious that a simple graph is assumed as a $\mathbb{T}$-gain graph with all positive gain 1's.

The \emph{positive inertia index} $i_{+}(\Phi)$, the \emph{negative inertia index} $i_{-}(\Phi)$ and the \emph{nullity} $\eta(\Phi)$ of $\Phi$ are defined to be the number of positive eigenvalues, negative eigenvalues and zero eigenvalues of $A(\Phi)$ including multiplicities, respectively. The \emph{rank} $r(\Phi)$ of an $n$-vertex $\mathbb{T}$-gain graph is defined as the rank of $A(\Phi)$. Obviously, $r(\Phi)=i_{+}(\Phi)+i_{-}(\Phi)=n-\eta(\Phi)$.

An \emph{induced subgraph} of $\Phi$ is an subgraph of $\Phi$ and each edge preserves the original gain in $\Phi$. For a vertex $v\in V(\Phi)$, we write $\Phi-v$ for the gain graph obtained from $\Phi$ by removing the vertex $v$ and all edges incident with $v$. For an induced subgraph $H$ of $\Phi$, let $\Phi-H$ be the subgraph obtained from $\Phi$ by deleting all vertices of $H$ and all incident edges. The \emph{degree} of a vertex $v$ for a gain graph $\Phi$ is the number of the vertices incident to $v$ in its underlying graph $G$. A vertex of a gain graph $\Phi=(G,\varphi)$ is called \emph{pendant} vertex if its degree is 1 in $\Phi$, and is called \emph{quasi-pendant} vertex if it is adjacent to a pendant vertex in $\Phi$. Denoted by $S_{n},~K_{n},~P_{n},~C_{n}$ a star, a complete graph, a path and a cycle all of order $n$, respectively. A graph is called \emph{trivial} if it has one vertex and no edges, it is sometimes denoted by $K_{1}$ or $P_{1}$.

A \emph{bicyclic} graph is a graph in which the number of edges equals the number of vertices plus one.
Let $G$ be a bicyclic graph, the \emph{base} of $G$ is the unique bicyclic subgraph of $G$ containing no pendant vertices.

Let $C_{p}$ and $C_{q}$ be two vertex-disjoint cycles and $v\in V(C_{p})$, $u\in V(C_{q})$, $P_{l}=v_{1}v_{2}\ldots v_{l}~(l\geq1)$ be a  path of length $l-1$. Let $\infty(p,l,q)$ (as shown in Fig. 1) be the graph obtained from $C_{p}$, $C_{q}$ and $P_{l}$ by identifying $v$ with $v_{1}$ and $u$ with $v_{l}$, respectively. The bicyclic graph containing $\infty(p,l,q)$ as its base is called an \emph{$\infty$-graph}. We denote $\infty^{\varphi}(p,l,q)$ be the $\mathbb{T}$-gain $\infty(p,l,q)$ graph, and the $\mathbb{T}$-gain bicyclic graph containing $\infty^{\varphi}(p,l,q)$ as its base is called a\emph{ $\mathbb{T}$-gain $\infty$-graph}.

Let $P_{p+2},P_{l+2},P_{q+2}$ be three paths, where min$\{p,l,q\}\geq0$ and at most one of $p,l,q$ is 0. Let $\theta(p,l,q)$ (as shown in Fig. 1) be the graph obtained from $P_{p+2}$, $P_{l+2}$ and $P_{q+2}$ by identifying the three initial vertices and terminal vertices. The bicyclic graph containing $\theta(p,l,q)$ as its base is called a \emph{$\theta$-graph}. We denote $\theta^{\varphi}(p,l,q)$ be the $\mathbb{T}$-gain $\theta(p,l,q)$ graph, and the $\mathbb{T}$-gain bicyclic graph containing $\theta^{\varphi}(p,l,q)$ as its base is called a\emph{ $\mathbb{T}$-gain $\theta$-graph}.

Recently there are many widely investigated research results about spectral-based graph invariants, such as inertia index of a graph \cite{YFW,YZF}, skew energy \cite{ACBR, LXLH} and skew-rank \cite{LXYG, LWZ, QHYG} of oriented graphs. Nathan Reff \cite{NR} defined the adjacency, incidence and Laplacian matrices of a complex unit gain graph. Some eigenvalue bounds for the adjacency and Laplacian matrices were present. Yu, Qu and Tu \cite{YQT} obtained some properties of inertia of a $\mathbb{T}$-gain graph. They characterized the $\mathbb{T}$-gain unicyclic graphs with small positive or negative index. Motivated by above, in this paper, we will investigate some characterizations about $\mathbb{T}$-gain bicyclic graphs.

The rest of this paper is organized as follows: in Section 2, we list some elementary lemmas and known results. In Section 3, we characterize the rank of $\mathbb{T}$-gain bicyclic graphs. In Section 4, we characterize the $\mathbb{T}$-gain bicyclic graphs with rank 2, 3 or 4.

\begin{figure}[htbp]
  \centering
  \includegraphics[scale=0.6]{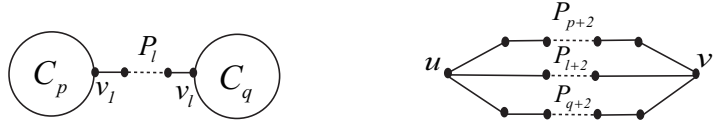}
\caption{$\infty(p,l,q)$ and $\theta(p,l,q)$.}
\end{figure}

\section{Preliminaries}
\noindent\begin{lemma}\label{le:2.1} \cite{YQT}
\begin{enumerate}[(a)]
  \item Let $\Phi=\Phi_{1}\bigcup \Phi_{2}\bigcup \ldots \bigcup \Phi_{t}$, where $\Phi_{1}, \Phi_{2}, \ldots, \Phi_{t}$ are connected components of a $\mathbb{T}$-gain graph $\Phi$. Then $i_{+}(\Phi)=\sum_{i=1}^{t}i_{+}(\Phi)_{i}$.
  \item Let $\Phi$ be a $\mathbb{T}$-gain graph on $n$ vertices. Then $i_{+}(\Phi)=0$ if and only if $\Phi$ is a graph without edges.
  \item Let $H^{\varphi}$ be an induced subgraph of a $\mathbb{T}$-gain graph $G^{\varphi}$. Then $i_{+}(H^{\varphi})\leq i_{+}(G^{\varphi}), i_{-}(H^{\varphi})\leq i_{-}(G^{\varphi})$.
\end{enumerate}
\end{lemma}

\noindent\begin{lemma}\label{le:2.2}\cite{YQT}
Let $\Phi=(G,\varphi)$ be a $\mathbb{T}$-gain graph containing a pendant vertex $v$ with the unique neighbor $u$. Then $i_{+}(\Phi)=i_{+}(\Phi-u-v)+1$, $i_{-}(\Phi)=i_{-}(\Phi-u-v)+1$, $i_{0}(\Phi)=i_{0}(\Phi-u-v)$. Moreover, $r(\Phi)=r(\Phi-u-v)+2$.
\end{lemma}

Let $C_{n}^{\varphi}$ be a weighted cycle with vertex set $\{v_{1}, v_{2},\ldots, v_{n}\}$ such that $v_{i}v_{i+1}\in E(C_{n}^{\varphi})$ ($1\leq i\leq n-1$), $v_{1}v_{n}\in E(C_{n}^{\varphi})$. Let $w_{i}=\varphi(v_{i}v_{i+1})$ and $w_{n}=\varphi(v_{n}v_{1})$.

In \cite{YQT}, Yu gave Definition 1 as follows about a $\mathbb{T}$-gain cycle. In fact, Definition 1 should be Definition 2.

\noindent\begin{definition}\label{de:1}\cite{YQT}
A $\mathbb{T}$-gain even cycle $C_{n}^{\varphi}$ is said to be of Type A if

$(-1)^{\frac{n}{2}}w_{n}=w_{1}w_{2}w_{3}\ldots w_{n-2}w_{n-1}$;

\noindent $C_{n}^{\varphi}$ is said to be of Type B if

$(-1)^{\frac{n}{2}}w_{n}\neq w_{1}w_{2}w_{3}\ldots w_{n-2}w_{n-1}$.

\noindent A $\mathbb{T}$-gain odd cycle $C_{n}^{\varphi}$ is said to be of Type C if

$Re\left((-1)^{\frac{n-1}{2}}w_{1}w_{2}w_{3}\ldots w_{n-2}w_{n-1}\overline{w_{n}}\right)>0$;

\noindent $C_{n}^{\varphi}$ is said to be of Type D if

$Re\left((-1)^{\frac{n-1}{2}}w_{1}w_{2}w_{3}\ldots w_{n-2}w_{n-1}\overline{w_{n}}\right)<0$;

\noindent $C_{n}^{\varphi}$ is said to be of Type E if

$Re\left((-1)^{\frac{n-1}{2}}w_{1}w_{2}w_{3}\ldots w_{n-2}w_{n-1}\overline{w_{n}}\right)=0$,

\end{definition}
where $Re(\cdot)$ is the real part of a complex number.

In the proof of Theorem 7 in \cite{YQT}, let $w_{i}=\varphi(v_{i}v_{i+1})$ $(1\leq i\leq n-1)$, $w_{n}=\varphi(v_{n}v_{1})$. The authors of \cite{YQT} gave the adjacency matrix $A(C_{n}^{\varphi})$ of $C_{n}^{\varphi}$ as follows:

\begin{displaymath}
A(C_{n}^{\varphi})=\left(
  \begin{array}{cccccccccccccc}
              0&            w_{1}&        0&           0&                \cdots&         0&    w_{n}\\
   \overline{w_{1}}&   0&     w_{2}&       0&                 \cdots&   0&   0\\
0&    \overline{ w_{2}}&    0&     w_{3}&              \cdots&   0&  0\\
 0&      0&     \overline{w_{3}}&          0&           \cdots&     0& 0\\
  \vdots&          \vdots&    \vdots&                \vdots&   \ddots&  \vdots& \vdots\\
            0&          0&   0&   0& \cdots&             0&            w_{n-1} \\
  \overline{ w_{n}}&  0&  0&  0&  \cdots&   \overline{w_{n-1}}&    0
  \end{array}
\right).
\end{displaymath}

From $A(C_{n}^{\varphi})$, we can see that $a_{1n}=w_{n}$ should be $a_{1n}=\overline{ w_{n}}$ and $a_{n1}=\overline{ w_{n}}$  should be $a_{n1}=w_{n}$. So, they made a mistake. In the following, we will give the right definition.

\noindent\begin{definition}\label{de:2}
A $\mathbb{T}$-gain even cycle $C_{n}^{\varphi}$ is said to be of Type A if

$(-1)^{\frac{n}{2}}\overline{w_{n}}=w_{1}w_{2}w_{3}\ldots w_{n-2}w_{n-1}$

(i.e., $\varphi(e_{n-1,n})+(-1)^{\frac{n-2}{2}}\frac{\varphi(e_{1,n})\varphi(e_{32})\varphi(e_{54})\cdots \varphi(e_{n-1,n-2})}{\varphi(e_{12})\varphi(e_{34})\cdots \varphi(e_{n-3,n-2})}=0$);

\noindent $C_{n}^{\varphi}$ is said to be of Type B if

$(-1)^{\frac{n}{2}}\overline{w_{n}}\neq w_{1}w_{2}w_{3}\ldots w_{n-2}w_{n-1}$

(i.e., $\varphi(e_{n-1,n})+(-1)^{\frac{n-2}{2}}\frac{\varphi(e_{1,n})\varphi(e_{32})\varphi(e_{54})\cdots \varphi(e_{n-1,n-2})}{\varphi(e_{12})\varphi(e_{34})\cdots \varphi(e_{n-3,n-2})}\neq0$).

\noindent A $\mathbb{T}$-gain odd cycle $C_{n}^{\varphi}$ is said to be of Type C if

$Re\left((-1)^{\frac{n-1}{2}}w_{1}w_{2}w_{3}\ldots w_{n-2}w_{n-1}w_{n}\right)>0$

(i.e., $Re\left((-1)^{\frac{n-1}{2}}\frac{\varphi(e_{32})\varphi(e_{54})\ldots \varphi(e_{n,n-1})\varphi(e_{1,n})}{\varphi(e_{12})\varphi(e_{34})\ldots \varphi(e_{n-2,n-1})}\right)>0$);

\noindent $C_{n}^{\varphi}$ is said to be of Type D if

$Re\left((-1)^{\frac{n-1}{2}}w_{1}w_{2}w_{3}\ldots w_{n-2}w_{n-1}w_{n}\right)<0$

(i.e., $Re\left((-1)^{\frac{n-1}{2}}\frac{\varphi(e_{32})\varphi(e_{54})\ldots \varphi(e_{n,n-1})\varphi(e_{1,n})}{\varphi(e_{12})\varphi(e_{34})\ldots \varphi(e_{n-2,n-1})}\right)<0$);

\noindent $C_{n}^{\varphi}$ is said to be of Type E if

$Re\left((-1)^{\frac{n-1}{2}}w_{1}w_{2}w_{3}\ldots w_{n-2}w_{n-1}w_{n}\right)=0$

(i.e., $Re\left((-1)^{\frac{n-1}{2}}\frac{\varphi(e_{32})\varphi(e_{54})\ldots \varphi(e_{n,n-1})\varphi(e_{1,n})}{\varphi(e_{12})\varphi(e_{34})\ldots \varphi(e_{n-2,n-1})}\right)=0$),

\end{definition}
where $Re(\cdot)$ is the real part of a complex number.

\noindent\begin{lemma}\label{le:2.3}\cite{YQT}
Let $C_{n}^{\varphi}$ be a $\mathbb{T}$-gain cycle of order $n$. Then

\begin{displaymath}
(i_{+}(C_{n}^{\varphi}),i_{-}(C_{n}^{\varphi}),i_{0}(C_{n}^{\varphi}))=\left\{\
        \begin{array}{ll}
          (\frac{n-2}{2},\frac{n-2}{2},2),& \emph{if}~C_{n}^{\varphi}\rm~is~of~Type~A,\\
          (\frac{n}{2},\frac{n}{2},0),& \emph{if}~C_{n}^{\varphi}\rm~is~of~Type~B,\\
          (\frac{n+1}{2},\frac{n-1}{2},0),& \emph{if}~C_{n}^{\varphi}\rm~is~of~Type~C,\\
          (\frac{n-1}{2},\frac{n+1}{2},0),& \emph{if}~C_{n}^{\varphi}\rm~is~of~Type~D,\\
          (\frac{n-1}{2},\frac{n-1}{2},1),& \emph{if}~C_{n}^{\varphi}\rm~is~of~Type~E.
        \end{array}
      \right.
\end{displaymath}

\end{lemma}

Two pendant vertices are called \emph{pendant twins} in a $\mathbb{T}$-gain graph $\Phi$ if they have the same neighbor in $\Phi$.

\noindent\begin{lemma}\label{le:2.4}
Let $u$, $v$ be pendant twins of a $\mathbb{T}$-gain graph $\Phi$. Then $r(\Phi)=r(\Phi-u)=r(\Phi-v)$.
\end{lemma}

\noindent\textbf{Proof.} Assume that all vertices in $\Phi$ are indexed by $\{v_{1},v_{2},\ldots,v_{n}\}$ with $v_{1}=v,v_{2}=u$. Without loss of generality, we assume $v, u$ are incident with $w=v_{3}$.
 Then it follows that

\begin{displaymath}
A(\Phi)=\left(
  \begin{array}{cccccccccccccc}
              0&            0&        \varphi(e_{13})&           0&                \cdots&         0 \\
              0&            0&        \varphi(e_{23})&           0&                 \cdots&   0\\
\varphi(e_{31})&    \varphi(e_{32})&         0&          \varphi(e_{34})&              \cdots&    \varphi(e_{3n}) \\
              0&            0&    \varphi(e_{43})&          0&           \cdots&       \varphi(e_{4n}) \\
  \vdots&          \vdots&    \vdots&                \vdots&   \ddots&  \vdots\\
            0&          0&    \varphi(e_{n3})&             \varphi(e_{n4})&                \cdots&   0 \\
  \end{array}
\right).
\end{displaymath}

By elementary row and column transformations on $A(\Phi)$, we have

\begin{displaymath}
r(A(\Phi))=r\left(
  \begin{array}{cccccccccccccc}
              0&             \varphi(e_{13})&           0&                \cdots&         0 \\
\varphi(e_{31})&             0&          \varphi(e_{34})&              \cdots&    \varphi(e_{3n})\\
              0&            \varphi(e_{43})&          0&           \cdots&       \varphi(e_{4n})\\
  \vdots&              \vdots&                \vdots&   \ddots&  \vdots\\
            0&              \varphi(e_{n3})&             \varphi(e_{n4})&                \cdots&   0\\
  \end{array}
\right).
\end{displaymath}

So, $r(\Phi)=r(\Phi-u)$. Similarly, $r(\Phi)=r(\Phi-v)$. \quad $\square$

From Lemmas \ref{le:2.1}, \ref{le:2.2} and \ref{le:2.4}, we can get the following two lemmas.

\noindent\begin{lemma}\label{le:2.5}

 Let $\Phi_{0}=(G_{0},\varphi_{0})$ be a $\mathbb{T}$-gain graph of order $n-p$ such that $u\in V(\Phi_{0})$. Let $\Phi_{1}=(G_{1},\varphi_{1})$ be the $\mathbb{T}$-gain graph obtained from $\Phi_{0}$ and $S_{p}^{\varphi_{3}}$ by inserting an edge between $u$ and the center $v_{1}$ of $S_{p}^{\varphi_{3}}$. Let
$\Phi_{2}=(G_{2},\varphi_{2})=\Phi_{1}-\{v_{1}v_{2},v_{1}v_{3},\ldots,v_{1}v_{p}\}+\{uv_{2},uv_{3},\ldots,uv_{p}\}$.
Then $r(\Phi_{1})\geq r(\Phi_{2})$.

\end{lemma}

\noindent\begin{lemma}\label{le:2.6}

Let $\Phi_{0}=(G_{0},\varphi_{0})$ be a $\mathbb{T}$-gain graph of order $n-l-t$ and $v_{1},v_{2}\in V(\Phi_{0})$. Assume that $\Phi_{1}=(G_{1},\varphi_{1})$ be the $\mathbb{T}$-gain graph obtained from $\Phi_{0}$, $S_{l+1}^{\varphi_{3}}$ and $S_{t+1}^{\varphi_{4}}$ by identifying $v_{1}$ with the center of $S_{l+1}^{\varphi_{3}}$, $v_{2}$ with the center of $S_{t+1}^{\varphi_{4}}$, respectively. Let $\Phi_{2}=(G_{2},\varphi_{2})$ be the $\mathbb{T}$-gain graph obtained from $\Phi_{0}$, $S_{l+t+1}^{\varphi_{5}}$ by identifying $v_{1}$ with the center of $S_{l+t+1}^{\varphi_{5}}$. Then $r(\Phi_{1})\geq r(\Phi_{2})$.

\end{lemma}

\noindent\begin{lemma}\label{le:2.7}

Let $\Phi_{1}=(G_{1},\varphi_{1})$ and $\Phi_{2}=(G_{2},\varphi_{2})$ be two $\mathbb{T}$-gain graphs and $v\in V(\Phi_{1})$, $u\in V(\Phi_{2})$. Let $P_{l}^{\varphi_{3}} (l\geq3)$ be a $\mathbb{T}$-gain path with two end-vertices $v_{1}$ and $v_{l}$. Let $\Phi=(G,\varphi)$ be the $\mathbb{T}$-gain graph obtained from $\Phi_{1}\bigcup \Phi_{2}\bigcup P_{l}^{\varphi_{3}}$ by identifying $v$ with $v_{1}$, $u$ with $v_{l}$, respectively. Let $\Phi'=(G',\varphi')$ be the  $\mathbb{T}$-gain graph obtained from $\Phi_{1}\bigcup \Phi_{2}$ by identifying $v$ with $u$ and adding $l-1$ pendant vertices to $v$. Then  $r(\Phi)\geq r(\Phi')$.

\end{lemma}

\noindent\textbf{Proof.} By Lemma \ref{le:2.2}, we have
$r(\Phi')=2+r(\Phi_{1}-v)+r(\Phi_{2}-u)$.

If $l=3$, by Lemmas \ref{le:2.1}(c) and \ref{le:2.2}, then we have $r(\Phi)\geq 2+r(\Phi_{1}-v)+r(\Phi_{2}-u)$.

If $l\geq4$, by Lemmas \ref{le:2.1}(c) and \ref{le:2.2}, then we have $r(\Phi)\geq 2+r(\Phi_{1}-v)+r(\Phi_{2})$. Note that $r(\Phi_{2})\geq r(\Phi_{2}-u)$, therefore $r(\Phi)\geq r(\Phi')$. \quad $\square$

\noindent\begin{theorem}\label{th:2.8}

Let $C_{n}^{\varphi_{0}}$ be a $\mathbb{T}$-gain cycle of order $n~(n\geq3)$ and $H=(G_{1},\varphi_{1})$ be a $\mathbb{T}$-gain graph of order $m~(m\geq1)$. Assume that $\Phi=(G,\varphi)$ is the $\mathbb{T}$-gain graph obtained by identifying a vertex of $C_{n}^{\varphi_{0}}$ with a vertex of $H$ (i.e., $V(C_{n}^{\varphi_{0}})\bigcap V(H)=v $). Let $F=(G_{2},\varphi_{2})$ be the induced subgraph obtained from $H$ by deleting the vertex $v$ and its incident edges. Then

\begin{displaymath}
\left\{\
        \begin{array}{ll}
          r(\Phi)=n-2+r(H),& \emph{if}~C_{n}^{\varphi}~\emph{is~of~Type~A},\\
          r(\Phi)=n+r(F), & \emph{if}~C_{n}^{\varphi}~\emph{is~of~Type~B},\\
          r(\Phi)=n-1+r(H),& \emph{if}~C_{n}^{\varphi}~\emph{is~of~Type~E},\\
          n-1+r(F)\leq r(\Phi)\leq n+r(H), & \emph{if}~C_{n}^{\varphi}~\emph{is~of~Type~C~or~D}.
        \end{array}
      \right.
\end{displaymath}
\end{theorem}

\noindent\textbf{Proof.} \textbf{Case 1.}
When $n$ is even, without loss of generality, we assume that $n=2p~(p\geq2)$ and $V(C_{n}^{\varphi_{0}})\bigcap V(H)=v_{2p}$. Then the adjacency matrix $A(\Phi)$ of $\Phi$ can be expressed as:
\begin{center}
\resizebox{13.3cm}{3.2cm}{%
$
\footnotesize A(\Phi)=\left(
  \begin{array}{cccccccccc}
    0 & \varphi(e_{12}) &  &  &  &  & \varphi(e_{1,2p})\\
    \varphi(e_{21}) & 0 & \varphi(e_{23}) &  \multicolumn{2}{l}{}{\Huge \textbf{0}}&  & 0\\
     & \varphi(e_{32}) & 0 & \varphi(e_{34}) &  &  & \vdots\\
     &  & \varphi(e_{43}) & \ddots & \ddots &  & \vdots &  \multicolumn{2}{l}{}{\Huge \textbf{0}}  \\
     \multicolumn{1}{c}{}{\Huge \textbf{0}}&  &  & \ddots & \ddots & \varphi(e_{2p-2,2p-1}) & 0 \\
     &  &  &  & \varphi(e_{2p-1,2p-2})& 0 & \varphi(e_{2p-1,2p}) &  &  &  \\
     \varphi(e_{2p,1})& 0 & \cdots & \cdots & 0 & \varphi(e_{2p,2p-1}) & 0 & \varphi(\alpha_{1}) & \cdots & \varphi(\alpha_{m-1}) \\
     &  &  &  &  &  & \overline{\varphi(\alpha_{1}}) &    \\
     &  &  \multicolumn{2}{c}{}{\Huge \textbf{0}}   &  &  &  \vdots &  \multicolumn{2}{l}{}{\Large \textbf{M}}  \\
     &  &  &  &  &  & \overline{\varphi(\alpha_{m-1}})\\
  \end{array}
\right),
$}
\end{center}
where $v_{i}v_{i+1}\in E(C_{n}^{\varphi_{0}})(1\leq i\leq 2p-1)$, $v_{1}v_{2p}\in E(C_{n}^{\varphi_{0}})$, $\alpha_{j}\in E(H),~j=1,2,\ldots,m-1,~M=A(F)$.

By elementary row and column transformations, it is obvious to show that $A(\Phi)$ is congruent to

\begin{displaymath}
A(\Phi_{1})=\left(
  \begin{array}{ccccccccccc}
        A_{1}&            &            &           &         &          &              &          &\\
             &      \ddots&            &           &         &          &              &          &\\
             &            &        A_{p-1}&           &         &          &              &          &\\
             &            &            &           &        0&    a&             0&    \cdots&              0&\\
             &            &            &           &  \overline{a}&         0&    \varphi(\alpha_{1})&    \cdots&     \varphi(\alpha_{m-1})\\
             &            &            &           &        0&              \overline{\varphi(\alpha_{1}})&          &               &\\
             &            &            &           &   \vdots&    \vdots&     \multicolumn{2}{l}{}{\Large \textbf{M}}\\
             &            &            &           &        0&               \overline{\varphi(\alpha_{m-1}})&         &         &     &\\
  \end{array}
\right),
\end{displaymath}

where
\begin{displaymath}
A_{i}=\left(
  \begin{array}{cccccccc}
        0&         \varphi(e_{2i-1,2i})\\
        \varphi(e_{2i,2i-1})&       0\\
  \end{array}
\right),~i=1,2,\ldots p-1,
\end{displaymath}

\begin{displaymath}
a=\varphi(e_{2p-1,2p})+(-1)^{\frac{2p-2}{2}}\frac{\varphi(e_{1,2p})\varphi(e_{32})\cdots \varphi(e_{2p-1,2p-2})}{\varphi(e_{12})\varphi(e_{34})\cdots \varphi(e_{2p-3,2p-2})}.
\end{displaymath}
So, we can get the following results

$(a)$ If $a=\varphi(e_{2p-1,2p})+(-1)^{\frac{2p-2}{2}}\frac{\varphi(e_{1,2p})\varphi(e_{32})\cdots \varphi(e_{2p-1,2p-2})}{\varphi(e_{12})\varphi(e_{34})\cdots \varphi(e_{2p-3,2p-2})}=0$, i.e.,
 $C_{n}^{\varphi}$ is of Type A, then we have
$r(\Phi)=2p-2+r(H)=n-2+r(H)$.

$(b)$ If $a=\varphi(e_{2p-1,2p})+(-1)^{\frac{2p-2}{2}}\frac{\varphi(e_{1,2p})\varphi(e_{32})\cdots \varphi(e_{2p-1,2p-2})}{\varphi(e_{12})\varphi(e_{34})\cdots \varphi(e_{2p-3,2p-2})}\neq 0$, i.e., $C_{n}^{\varphi}$ is of Type B, then we have
$r(\Phi)=2p+r(F)=n+r(F)$.

\textbf{Case 2.} When $n$ is odd, without loss of generality, we assume that $n=2p-1~(p\geq2)$ and $V(C_{n}^{\varphi_{0}})\bigcap V(H)=v_{2p-1}$. By elementary row and column transformations, we can show $A(\Phi)$ is congruent to

\begin{displaymath}
A(\Phi_{1})=\left(
  \begin{array}{cccccc}
        A_{1}&         &         &                &     &\\
         &        A_{2}&         &                &           &\\
         &             &   \ddots&                &           &\\
         &             &         &           A_{p-1}&     &\\
         &             &         &                &          C&\\
  \end{array}
\right),
\end{displaymath}
where
\begin{displaymath}
  A_{i}=\left(
  \begin{array}{cc}
        0&        \varphi(e_{2i-1,2i})       \\
          \varphi(e_{2i,2i-1})&       0     \\
 \end{array}
  \right),
  \end{displaymath}
 \begin{displaymath}
  C=\left(
  \begin{array}{ccccccccccccc}
      a&    \varphi(\alpha_{1})&   \cdots&   \varphi(\alpha_{m-1})\\
      \overline{\varphi(\alpha_{1})}         \\
   \vdots&        \multicolumn{2}{l}{}{\Large \textbf{M}}\\
    \overline{\varphi(\alpha_{m-1})}         \\
  \end{array}
\right),
\end{displaymath}

$a=2Re\left((-1)^{\frac{2p-2}{2}}\frac{\varphi(e_{1,2p-1})\varphi(e_{32})\cdots \varphi(e_{2p-1,2p-2})}{\varphi(e_{12})\varphi(e_{34})\cdots \varphi(e_{2p-3,2p-2})}\right)$,
for $i=1,2,\ldots, p-1$. So, we can get the following results

$(c)$ If $a=0$, i.e., $C_{n}^{\varphi}$ is of Type E, then we have
$r(\Phi)=2p-2+r(H)=n-1+r(H)$.

$(d)$ If $a\neq0$, i.e., $C_{n}^{\varphi}$ is of Type C or D, by Lemmas \ref{le:2.1}(c), \ref{le:2.2} and $r(C)\leq r(H)+1$, then we have
$n-1+r(F)\leq r(\Phi)\leq n+r(H)$. \quad $\square$

\section{The rank of $\mathbb{T}$-gain bicyclic graphs}

In this section, we shall characterize the rank of $\mathbb{T}$-gain bicyclic graphs. First, we will give two theorems.

\begin{figure}[htbp]
  \centering
  \includegraphics[scale=0.6]{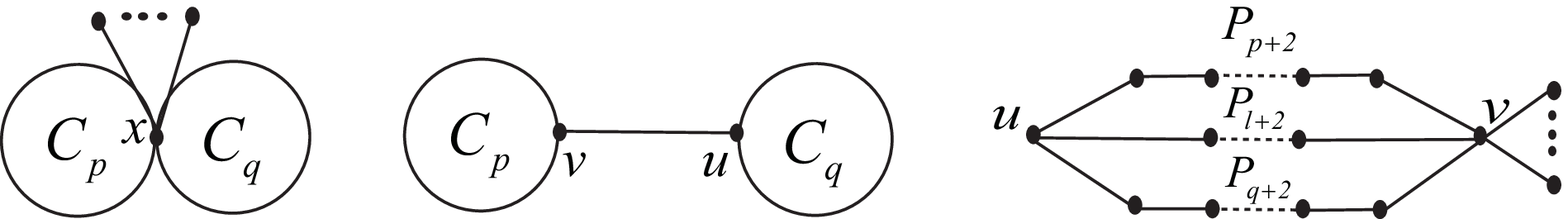}
\caption{$G^{\ast}$,~$\infty(p,2,q)$ and $G^{\ast\ast}$.}
\end{figure}

\noindent\begin{theorem}\label{th:3.1}
Let $\Phi=(G,\varphi)$ be a $\mathbb{T}$-gain bicyclic graph of order $n~(n\geq p+q)$ with pendant vertices containing two edge disjoint cycles $C_{p}^{\varphi_{1}}$ and $C_{q}^{\varphi_{2}}$. Then
\begin{displaymath}
r(\Phi)\geq\left\{\
        \begin{array}{ll}
          p+q,&\emph{if}~p,~q~\emph{are~odd},\\
          p+q-2,&\emph{if}~p,~q~\emph{are~even},\\
          p+q-1,&\emph{otherwise}.
        \end{array}
      \right.
\end{displaymath}

\noindent This bound is sharp.
\end{theorem}

\noindent\textbf{Proof.} Let $\Phi^{\ast}=(G^{\ast},\varphi^{\ast})$, where $G^{\ast}$ (as shown in Fig. 2) is the underlying graph of $\Phi^{\ast}$. By Lemma \ref{le:2.2}, we have $r(\Phi^{\ast})=2+r(P_{p-1}^{\varphi_{1}})+r(P_{q-1}^{\varphi_{2}})$.

Let $\Phi=(G,\varphi)$ be a $\mathbb{T}$-gain bicyclic graph of order $n~(n\geq p+q)$ with pendant vertices containing two edge disjoint cycles $C_{p}^{\varphi_{1}}$ and $C_{q}^{\varphi_{2}}$. In the following we shall prove that $r(\Phi)\geq r(\Phi^{\ast})$. By Lemmas \ref{le:2.5}, \ref{le:2.6} and \ref{le:2.7}, we will distinguish the following two cases.

\textbf{Case 1.} $\Phi$ is a $\mathbb{T}$-gain bicyclic graph obtained from $\infty^{\varphi}(p,1,q)$ by attaching $n-p-q+1~(n\geq p+q)$ pendant vertices to a vertex, different from $x$, of $C_{p}^{\varphi_{1}}$ or $C_{q}^{\varphi_{2}}$, where $x$ is the common vertex of $C_{p}^{\varphi_{1}}$ and $C_{q}^{\varphi_{2}}$.

Without loss of generality, we suppose all the pendant vertices are attached at $C_{p}^{\varphi_{1}}$. Then we have

\begin{displaymath}
r(\Phi)=2+r(P_{p-1}^{\varphi_{1}})+\left\{\
        \begin{array}{ll}
          r(P_{q-1}^{\varphi_{2}}),&\textrm{if}~p~\textrm{is~odd},\\
         r(P_{q-1}^{\varphi_{2}})\textrm{~or}~r(C_{q}^{\varphi_{2}}),&\textrm{if}~p~\textrm{is~even}.

        \end{array}
      \right.
\end{displaymath}

From Lemma \ref{le:2.1}(c) and Lemma \ref{le:2.3}, we have $r(C_{q}^{\varphi_{2}})\geq r(P_{q-1}^{\varphi_{2}})$. So we have $r(\Phi)\geq r(\Phi^{\ast})$.

\textbf{Case 2.} $\Phi$ is a $\mathbb{T}$-gain bicyclic graph obtained from $\infty^{\varphi}(p,2,q)$ by attaching $n-p-q~(n\geq p+q+1)$ pendant vertices to a vertex $w$ of $C_{p}^{\varphi_{1}}$ or $C_{q}^{\varphi_{2}}$.

Without loss of generality, we suppose all the pendant vertices are attached at $C_{p}^{\varphi_{1}}$. If $w=v$ (where $v$ is shown in Fig.~2), then we have $r(\Phi)=2+r(P_{p-1}^{\varphi_{1}})+r(C_{q}^{\varphi_{2}})\geq r(\Phi^{\ast})$. If $w\neq v$, then we have

\begin{displaymath}
r(\Phi)=2+\left\{\
        \begin{array}{ll}
          r(P_{p-1}^{\varphi_{1}})+r(C_{q}^{\varphi_{2}}),&\textrm{if}~p~\textrm{is~odd},\\
         p+r(P_{q-1}^{\varphi_{2}})~\textrm{or}
         ~r(P_{p-1}^{\varphi_{1}})+r(C_{q}^{\varphi_{2}}),&\textrm{if}~p~\textrm{is~even}.
        \end{array}
      \right.
\end{displaymath}

So we have $r(\Phi)\geq r(\Phi^{\ast})$. \quad $\square$

Let $u,v$ be two vertices in $\theta^{\varphi}(p,l,q)$ (as show in Fig.~1). Let $\Phi^{\ast\ast}=(G^{\ast\ast},\varphi)$ be the $\mathbb{T}$-gain bicyclic graph with $n-p-q-l-2~(n\geq p+q+l+3)$ pendant vertices attached to $v$ in  $\theta^{\varphi}(p,l,q)$ (as show in Fig. 2).

\noindent\begin{theorem}\label{th:3.2}

Let $\Phi$ be a $\mathbb{T}$-gain bicyclic graph of order $n~(n\geq p+q+l+3)$ with pendant vertices containing $\theta^{\varphi}(p,l,q)$ as its base. If $plq\neq0$, then

\begin{displaymath}
r(\Phi)\geq\left\{
        \begin{array}{ll}
          p+q+l+2,~ &\emph{if}~p+q+l~\emph{is~even}, \\
          p+q+l+1,&\emph{if}~p,~q,~l~\emph{are~odd,~or}~p~\emph{is~odd},~C_{q+l+2}^{\varphi}~\emph{is}\\
          &\emph{of~Type~A},\\
          p+q+l+3, &\rm otherwise. \\
        \end{array}
      \right.
\end{displaymath}

\noindent This bound is sharp.
\end{theorem}

\noindent\textbf{Proof.} Let $\Phi^{\ast\ast}=\Phi_{1}$ be a $\mathbb{T}$-gain graph with $G^{\ast\ast}$ as the underlying graph. By Lemma \ref{le:2.2}, we have

\begin{displaymath}
r(\Phi_{1})=\left\{\
        \begin{array}{ll}
          2+p+r(P_{q+l+1}^{\varphi_{0}}),&\textrm{if}~p~\textrm{is~even},\\
         3+p+r(P_{q}^{\varphi_{1}})+r(P_{l}^{\varphi_{2}}),&\textrm{if}~p~\textrm{is~odd}.
        \end{array}
      \right.
\end{displaymath}

Let $\Phi_{2}$ be a $\mathbb{T}$-gain bicyclic graph of order $n~(n\geq p+q+l+3)$ with pendant vertices containing $\theta^{\varphi}(p,l,q)$ as its base.  Consider the graph $\Phi_{2}$ in which all $n-p-q-l-2$ pendant vertices are attached to a vertex, different from $u$ and $v$ of $\theta^{\varphi}(p,l,q)$.
Without loss of generality, assume that $n-p-q-l-2$ pendant vertices are attached to a vertex of $P_{p+2}^{\varphi_{3}}$ in $\Phi_{2}$. By Lemma \ref{le:2.2}, we have

\begin{displaymath}
r(\Phi_{2})=\left\{\
        \begin{array}{ll}
          2+p+r(P_{q+l+1}^{\varphi_{0}}),&\textrm{if}~p~\textrm{is~even},\\
         3+p+r(P_{q}^{\varphi_{1}})+r(P_{l}^{\varphi_{2}})~\textrm{or}~1+p+r(C_{q+l+2}^{\varphi_{4}}),&\textrm{if}~p~\textrm{is~odd}.
        \end{array}
      \right.
\end{displaymath}

Combing the ranks of $\Phi_{1}$, $\Phi_{2}$ and Lemmas \ref{le:2.2} and \ref{le:2.3}, we can get the results.\quad $\square$

Next we consider the special case that one of $p,q,l$ is zero. Without loss of generality, we may assume that $l=0$. By a similar discussion as in the proof of Theorem \ref{th:3.2}, we can get the following result.

\noindent\begin{theorem}\label{th:3.3}

Let $\Phi$ be a $\mathbb{T}$-gain bicyclic graph of order $n$ with pendant vertices containing $\theta^{\varphi}(p,0,q)$ as its base. Then

\begin{displaymath}
r(\Phi)\geq\left\{
        \begin{array}{ll}
          2+p+q,~ &\emph{if}~p+q~\emph{is~even}, \\
          1+p+q, &\rm otherwise.
        \end{array}
      \right.
\end{displaymath}
\noindent This bound is sharp.
\end{theorem}

Combining Theorems \ref{th:3.1}, \ref{th:3.2} and \ref{th:3.3}, we can get the following theorem.

\noindent\begin{theorem}\label{th:3.4}

Let $\Phi=(G,\varphi)$ be a $\mathbb{T}$-gain bicyclic graph of order $n$ with pendant vertices. Then

\begin{enumerate}[(1)]
  \item If $\Phi$ is a $\mathbb{T}$-gain $\infty$-graph, then $r(\Phi)\geq6$.
  \item If $\Phi$ is a $\mathbb{T}$-gain $\theta$-graph, then $r(\Phi)\geq4$.
\end{enumerate}

\end{theorem}

\begin{figure}[htbp]
\begin{minipage}[hbt]{1\columnwidth}
\textbf{Table 1}

 {\small { The gain conditions for each gain graph in Theorem \ref{th:4.1} with rank 2, 3 or 4.}}
  \vskip1mm
 {\footnotesize
\begin{tabular}{p{2cm}p{11cm}}
\hline
Rank of $G^{\varphi}$& Gain graph $G^{\varphi}$ and its gain conditions \\
\hline
$r(G^{\varphi})=2$ & $G_{5}^{\varphi}$, the subgraph induced on vertices $1,~2,~3$ is of Type E and the subgraph induced on vertices $1,~2,~4,~3$ is of Type A.

$G_{9}^{\varphi}$, the subgraphs induced on vertices $1,~2,~3,~4$ and $1,~2,~5,~4$ are of Type A.

\\

$r(G^{\varphi})=3$ & $G_{5}^{\varphi}$, the subgraph induced on vertices $1,~2,~3$ is of Type C or D, and the subgraph induced on vertices $1,~2,~4,~3$ is of Type A.

\\

$r(G^{\varphi})=4$ & $G_{1}^{\varphi}$,  $Re\left(-\varphi(e_{15})\frac{\varphi(e_{52})}{\varphi(e_{12})}\right)+Re\left(-\varphi(e_{35})\frac{\varphi(e_{54})}{\varphi(e_{34})}\right)=0$.

$G_{2}^{\varphi}$, the subgraph induced on vertices $3,~4,~5,~6$ is of Type A and the subgraph induced on vertices $1,~2,~3$ is of Type E.

$G_{3}^{\varphi}$, the subgraphs induced on vertices $1,~2,~4,~3$ and $4,~5,~6,~7$ are of Type A.

$G_{5}^{\varphi}$, the subgraph induced on vertices $1,~2,~3$ is of Type C, or D, or E and the subgraph induced on vertices $1,~2,~4,~3$ is of Type B.

$G_{6}^{\varphi}$,  $Re\left(-\varphi(e_{15})\frac{\varphi(e_{52})}{\varphi(e_{12})}\right)+Re\left(\frac{\varphi(e_{15})\varphi(e_{32})\varphi(e_{54})}{\varphi(e_{12})\varphi(e_{34})}\right)=0$.

$G_{7}^{\varphi}$, the subgraph induced on vertices $1,~2,~6$ is of Type E and the subgraph induced on vertices $1,~2,~3,~4,~5,~6$ is of Type A.

$G_{8}^{\varphi}$,
$\varphi(e_{16})\varphi(e_{32})\varphi(e_{54})-$ $\varphi(e_{16})\varphi(e_{34})\varphi(e_{52})+$$\varphi(e_{12})\varphi(e_{34})\varphi(e_{56})=0$.

$G_{9}^{\varphi}$, the subgraph induced on vertices $1,~2,~3,~4$ is of Type A and the subgraph induced on vertices $1,~2,~5,~4$ is of Type B, or the subgraph induced on vertices $1,~2,~3,~4$ is of Type B and the subgraph induced on vertices $1,~2,~5,~4$ is of Type A or B.

$G_{10}^{\varphi}$, the subgraph induced on vertices $1,~2,~3,~4,~5$ is of Type E and the subgraph induced on vertices $1,~2,~6,~5$ is of Type A.

$G_{11}^{\varphi}$, the subgraphs induced on vertices $1,~2,~3,~4,~5,~6$ and $1,~2,~7,~6$ are of Type A.

\\

\hline
\end{tabular} }
\end{minipage}
\end{figure}

\begin{figure}[htbp]
  \centering
  \includegraphics[scale=0.6]{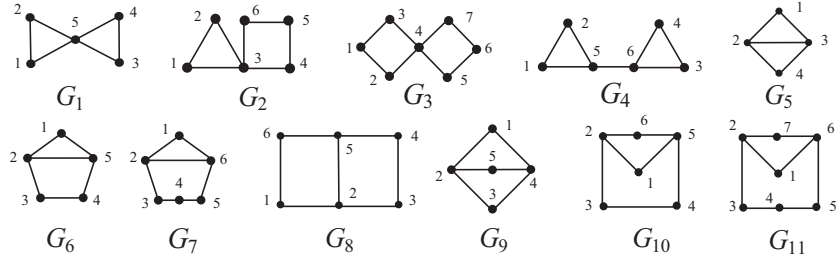}
\caption{The eleven graphs in Theorem 4.1.}
\end{figure}
\section{The $\mathbb{T}$-gain bicyclic graphs with rank 2, 3 or 4}

In this section, we shall characterize the $\mathbb{T}$-gain bicyclic graphs with rank 2, 3 or 4.

\noindent\begin{theorem}\label{th:4.1}

Let $\Phi=(G,\varphi)$ be a $\mathbb{T}$-gain bicyclic graph without pendant vertex. Then

\begin{enumerate}[(1)]
  \item $r(\Phi)=2$ if and only if $\Phi$ with one of $G_{i}$'s $(i=5,~9)$ (as shown in Fig. 3)  as its underlying graph, and the corresponding gain conditions are as shown in Table 1.
  \item $r(\Phi)=3$ if and only if $\Phi$ with $G_{5}$ (as shown in Fig. 3)  as its underlying graph, and the corresponding gain conditions are as shown in Table 1.
  \item $r(\Phi)=4$ if and only if $\Phi$ with one of $G_{i}$'s $(i=1-3,~5-11)$ (as shown in Fig. 3)  as its underlying graph, and the corresponding gain conditions are as shown in Table 1.
\end{enumerate}

\end{theorem}

\noindent\textbf{Proof.} Let $\Phi=(G,\varphi)$ be a $\mathbb{T}$-gain bicyclic graph without pendant vertices. Then $G$ is a underlying graph of $\Phi$, and $G^{\varphi}$ is a base.

\textbf{Case 1.} $G$ is an $\infty$-graph $\infty(p,l,q)$.

Without loss of generality, we suppose that $p\leq q$. If $(p,l,q)\in$ $\{(3,1,3),$ $~(3,1,4),~(4,1,4),~(3,2,3)\}$, then $G$ is one of the graphs $G_{1}-G_{{4}}$ in Fig. 3.

\begin{displaymath}
A(G_{1}^{\varphi})=\left(
  \begin{array}{cccccccccccccc}
            0&        \varphi(e_{12})&        0&           0&                    \varphi(e_{15})\\
     \varphi(e_{21})&               0&        0&           0&                    \varphi(e_{25})\\
            0&                      0&        0&        \varphi(e_{34})&         \varphi(e_{35})\\
            0&                      0&  \varphi(e_{43})&   0&                    \varphi(e_{45})\\
          \varphi(e_{51})&  \varphi(e_{52})&     \varphi(e_{53})&     \varphi(e_{54})&         0\\
  \end{array}
\right).
\end{displaymath}
By simple calculations, we have
\begin{displaymath}
r(G_{1}^{\varphi})=r\left(
  \begin{array}{cccccccccccccc}
            0&         \varphi(e_{12}) \\
            \varphi(e_{21})&          0 \\
  \end{array}
\right)+r\left(
  \begin{array}{cccccccccccccc}
            0&         \varphi(e_{34}) \\
            \varphi(e_{43})&          0 \\
  \end{array}
\right)+r(a+\overline{a}+b+\overline{b}),
\end{displaymath}

where $a=-\varphi(e_{15})\frac{\varphi(e_{52})}{\varphi(e_{12})}$, $b=-\varphi(e_{35})\frac{\varphi(e_{54})}{\varphi(e_{34})}$. So $r(G_{1}^{\varphi})=4$ if and only if $a+\overline{a}+b+\overline{b}=0$, i.e., $Re\left(-\varphi(e_{15})\frac{\varphi(e_{52})}{\varphi(e_{12})}\right)+Re\left(-\varphi(e_{35})\frac{\varphi(e_{54})}{\varphi(e_{34})}\right)=0$.

By Lemma \ref{le:2.3} and Theorem \ref{th:2.8}, we have $r(G_{2}^{\varphi})=4$ if and only if the subgraph induced on vertices $3,~4,~5,~6$ is of Type A and the subgraph induced on vertices $1,~2,~3$ is of Type E.

By Lemma \ref{le:2.3} and Theorem \ref{th:2.8}, we have $r(G_{3}^{\varphi})=4$ if and only if all the subgraphs induced on vertices $1,~2,~4,~3$ and $4,~5,~6,~7$ are of Type A.

\begin{displaymath}
A(G_{4}^{\varphi})=\left(
  \begin{array}{cccccccccccccc}
            0&        \varphi(e_{12})&        0&           0&                    \varphi(e_{15})&    0\\
     \varphi(e_{21})&               0&        0&           0&                    \varphi(e_{25})&    0\\
            0&                      0&        0&        \varphi(e_{34})&     0&    \varphi(e_{36})\\
            0&                      0&  \varphi(e_{43})&   0&            0&        \varphi(e_{46})\\
          \varphi(e_{51})&  \varphi(e_{52})&     0&    0&  0&                      \varphi(e_{56})\\
          0&  0&  \varphi(e_{63})&    \varphi(e_{64})&    \varphi(e_{65})&   0
  \end{array}
\right).
\end{displaymath}
By simple calculations, we have

\begin{displaymath}
r(G_{4}^{\varphi})=r\left(
  \begin{array}{cccccccccccccc}
            0&         \varphi(e_{12}) \\
            \varphi(e_{21})&          0 \\
  \end{array}
\right)+r\left(
  \begin{array}{cccccccccccccc}
            0&         \varphi(e_{34}) \\
            \varphi(e_{43})&          0 \\
  \end{array}
\right)+r\left(
  \begin{array}{cccccccccccccc}
            a+\overline{a}&         \varphi(e_{56}) \\
            \varphi(e_{65})&         b+\overline{b} \\
  \end{array}
  \right),
\end{displaymath}
where $a=-\varphi(e_{15})\frac{\varphi(e_{52})}{\varphi(e_{12})}$, $b=-\varphi(e_{36})\frac{\varphi(e_{64})}{\varphi(e_{34})}$.

So, $5\leq r(G_{4}^{\varphi})\leq6$ since $\varphi(e_{56})\neq0$ and $\varphi(e_{65})\neq0$.

If $p\geq3,~l\geq1,~q\geq5$ or $p\geq3,~l\geq2,~q\geq4$, by Lemmas \ref{le:2.1}, \ref{le:2.2} and \ref{le:2.3}, then we have $r(\Phi)\geq6$.

\textbf{Case 2.} $G$ is a $\theta$-graph $\theta(p,l,q)$.

Without loss of generality, we suppose that $p\leq l\leq q$. If $(p,l,q)\in$ $\{(0,1,1)~(0,1,2),~(0,1,3),$ $~(0,2,2),~(1,1,1),~(1,1,2),~(1,1,3)\}$, then $G$ is one of the graphs $G_{5}-G_{11}$ in Fig. 3.

\begin{displaymath}
A(G_{5}^{\varphi})=\left(
  \begin{array}{cccccccccccccc}
            0&        \varphi(e_{12})&        \varphi(e_{13})&           0\\
     \varphi(e_{21})&               0&        \varphi(e_{23})&         \varphi(e_{24})\\
            \varphi(e_{31})&                      \varphi(e_{32})&        0&       \varphi(e_{34})\\
            0&        \varphi(e_{42})&  \varphi(e_{43})&   0\\
  \end{array}
\right).
\end{displaymath}
By simple calculations, we have
\begin{displaymath}
r(G_{5}^{\varphi})=r\left(
  \begin{array}{cccccccccccccc}
            0&         \varphi(e_{12}) \\
            \varphi(e_{21})&          0 \\
  \end{array}
\right)+r\left(
  \begin{array}{cccccccccccccc}
            a+\overline{a}&         \overline{b} \\
                         b&                    0 \\
  \end{array}
\right),
\end{displaymath}

where $a=-\varphi(e_{13})\frac{\varphi(e_{32})}{\varphi(e_{12})}$, $b=\varphi(e_{43})-\varphi(e_{13})\frac{\varphi(e_{42})}{\varphi(e_{12})}$.

So, we have

\begin{displaymath}
r(G_{5}^{\varphi})=\left\{
        \begin{array}{ll}
          2,~\textrm{if}~Re(a)=0,~b=0,\\
          3,~\textrm{if}~Re(a)\neq0,~b=0,\\
          4,~\textrm{if}~Re(a)=0,~b\neq0~\textrm{or}~Re(a)\neq0,~b\neq0.
        \end{array}
      \right.
\end{displaymath}

i.e., $r(G_{5}^{\varphi})=2$ if and only if  the subgraph induced on vertices $1,~2,~3$ is of Type E and the subgraph induced on vertices $1,~2,~4,~3$ is of Type A.

 $r(G_{5}^{\varphi})=3$ if and only if  the subgraph induced on vertices $1,~2,~3$ is of Type C or D and the subgraph induced on vertices $1,~2,~4,~3$ is of Type A.

 $r(G_{5}^{\varphi})=4$ if and only if the subgraph induced on vertices $1,~2,~3$ is of Type E and the subgraph induced on vertices $1,~2,~4,~3$ is of Type B, or the subgraph induced on vertices $1,~2,~3$ is of Type C or D and the subgraph induced on vertices $1,~2,~4,~3$ is of Type B.

\begin{displaymath}
A(G_{6}^{\varphi})=\left(
  \begin{array}{cccccccccccccc}
            0&        \varphi(e_{12})&       0&  0& \varphi(e_{15})\\
     \varphi(e_{21})&               0&        \varphi(e_{23})&    0&     \varphi(e_{25})\\
          0&  \varphi(e_{32})&               0&       \varphi(e_{34})&         0\\
            0&     0&   \varphi(e_{43})&  0&  \varphi(e_{45})\\
            \varphi(e_{51})& \varphi(e_{52})& 0&  \varphi(e_{54})&  0
  \end{array}
\right).
\end{displaymath}
By simple calculations, we have
\begin{displaymath}
r(G_{6}^{\varphi})=r\left(
  \begin{array}{cccccccccccccc}
            0&         \varphi(e_{12}) \\
            \varphi(e_{21})&          0 \\
  \end{array}
\right)+r\left(
  \begin{array}{cccccccccccccc}
           0&         \varphi(e_{34}) \\
           \varphi(e_{43})&                    0 \\
  \end{array}
\right)+r(a+\overline{a}+b+\overline{b}),
\end{displaymath}

where $a=-\varphi(e_{15})\frac{\varphi(e_{52})}{\varphi(e_{12})}$, $b=\frac{\varphi(e_{15})\varphi(e_{32})\varphi(e_{54})}{\varphi(e_{12})\varphi(e_{34})}$.

So, $r(G_{6}^{\varphi})=4$ if and only if $Re\left(-\varphi(e_{15})\frac{\varphi(e_{52})}{\varphi(e_{12})}\right)+Re\left(\frac{\varphi(e_{15})\varphi(e_{32})\varphi(e_{54})}{\varphi(e_{12})\varphi(e_{34})}\right)=0$.

\begin{displaymath}
A(G_{7}^{\varphi})=\left(
  \begin{array}{cccccccccccccc}
            0&        \varphi(e_{12})&    0&   0&  0& \varphi(e_{16})\\
     \varphi(e_{21})& 0& \varphi(e_{23})& 0&  0& \varphi(e_{26})\\
          0&  \varphi(e_{32})&               0&       \varphi(e_{34})&   0&      0\\
            0&     0&   \varphi(e_{43})&  0&  \varphi(e_{45})&  0\\
            0&  0&  0&  \varphi(e_{54})&  0&  \varphi(e_{56})\\
            \varphi(e_{61})& \varphi(e_{62})& 0&  0& \varphi(e_{65})&  0
  \end{array}
\right).
\end{displaymath}
By simple calculations, we have
\begin{displaymath}
r(G_{7}^{\varphi})=r\left(
  \begin{array}{cccccccccccccc}
            0&         \varphi(e_{12}) \\
            \varphi(e_{21})&          0 \\
  \end{array}
\right)+r\left(
  \begin{array}{cccccccccccccc}
           0&         \varphi(e_{34}) \\
           \varphi(e_{43})&                    0 \\
  \end{array}
\right)+r\left(
  \begin{array}{cccccccccccccc}
           0&         b \\
           \overline{b}&                   a+\overline{a} \\
  \end{array}
\right),
\end{displaymath}

where $a=-\varphi(e_{16})\frac{\varphi(e_{62})}{\varphi(e_{12})}$, $b=\varphi(e_{56})+\frac{\varphi(e_{16})\varphi(e_{32})\varphi(e_{54})}{\varphi(e_{12})\varphi(e_{34})}$.

So, $r(G_{7}^{\varphi})=4$ if and only if the subgraph induced on vertices $1,~2,~6$ is of Type E and the subgraph induced on vertices $1,~2,~3,~4,~5,~6$ is of Type A.

\begin{displaymath}
A(G_{8}^{\varphi})=\left(
  \begin{array}{cccccccccccccc}
            0&        \varphi(e_{12})&    0&   0&  0& \varphi(e_{16})\\
     \varphi(e_{21})& 0& \varphi(e_{23})& 0& \varphi(e_{25})&  0\\
          0&  \varphi(e_{32})&               0&       \varphi(e_{34})&   0&      0\\
            0&     0&   \varphi(e_{43})&  0&  \varphi(e_{45})&  0\\
            0&  \varphi(e_{52})&  0&  \varphi(e_{54})&  0&  \varphi(e_{56})\\
            \varphi(e_{61})& 0& 0&  0& \varphi(e_{65})&  0
  \end{array}
\right).
\end{displaymath}
By simple calculations, we have

\begin{displaymath}
r(G_{8}^{\varphi})=r\left(
  \begin{array}{cccccccccccccc}
            0&         \varphi(e_{12}) \\
            \varphi(e_{21})&          0 \\
  \end{array}
\right)+r\left(
  \begin{array}{cccccccccccccc}
           0&         \varphi(e_{34}) \\
           \varphi(e_{43})&                    0 \\
  \end{array}
\right)+r\left(
  \begin{array}{cccccccccccccc}
           0&         a \\
           \overline{a}&                  0 \\
  \end{array}
\right),
\end{displaymath}
where $a=\frac{\varphi(e_{16})\varphi(e_{32})\varphi(e_{54})-\varphi(e_{16})\varphi(e_{34})\varphi(e_{52})+\varphi(e_{12})\varphi(e_{34})\varphi(e_{56})}{\varphi(e_{12})\varphi(e_{34})}$.

So, $r(G_{8}^{\varphi})=4$ if and only if $a=0$, i.e.,

$\varphi(e_{16})\varphi(e_{32})\varphi(e_{54})-\varphi(e_{16})\varphi(e_{34})\varphi(e_{52})+\varphi(e_{12})\varphi(e_{34})\varphi(e_{56})=0$.

\begin{displaymath}
A(G_{9}^{\varphi})=\left(
  \begin{array}{cccccccccccccc}
            0&        \varphi(e_{12})&      0&  \varphi(e_{14})&           0\\
     \varphi(e_{21})&               0&        \varphi(e_{23})&     0&    \varphi(e_{25})\\
         0&   \varphi(e_{32})&        0&              \varphi(e_{34})&        0\\
          \varphi(e_{41})&   0&        \varphi(e_{43})&  0&    \varphi(e_{45})\\
          0&   \varphi(e_{52})&  0&    \varphi(e_{54})&  0
  \end{array}
\right).
\end{displaymath}
By simple calculations, we have
\begin{displaymath}
r(G_{9}^{\varphi})=r\left(
  \begin{array}{cccccccccccccc}
            0&         \varphi(e_{12}) \\
            \varphi(e_{21})&          0 \\
  \end{array}
\right)+r\left(
  \begin{array}{cccccccccccccc}
               0&         a&                  0 \\
    \overline{a}&         0&        \overline{b}\\
               0&         b&                   0
  \end{array}
\right),
\end{displaymath}

where $a=\varphi(e_{34})-\varphi(e_{14})\frac{\varphi(e_{32})}{\varphi(e_{12})}$, $b=\varphi(e_{54})-\varphi(e_{14})\frac{\varphi(e_{52})}{\varphi(e_{12})}$.

So, we have

\begin{displaymath}
r(G_{9}^{\varphi})=\left\{
        \begin{array}{ll}
          2,~\textrm{if}~a=0,~b=0,\\
          4,~\textrm{if}~a=0,~b\neq0,~\textrm{or}~a\neq0,~b=0,~\textrm{or}~a\neq0,~b\neq0.
        \end{array}
      \right.
\end{displaymath}

i.e., $r(G_{9}^{\varphi})=2$ if and only if  the subgraphs induced on vertices $1,~2,~3,~4$ and  $1,~2,~5,~4$ are of Type A.

$r(G_{9}^{\varphi})=4$ if and only if  the subgraph induced on vertices $1,~2,~3,~4$ is of Type A and the subgraph induced on vertices $1,~2,~5,~4$ is of Type B, or the subgraph induced on vertices $1,~2,~3,~4$ is of Type B and the subgraph induced on vertices $1,~2,~5,~4$ is of Type A or B.

\begin{displaymath}
A(G_{10}^{\varphi})=\left(
  \begin{array}{cccccccccccccc}
            0&        \varphi(e_{12})&    0&   0&  \varphi(e_{15})&    0\\
     \varphi(e_{21})& 0& \varphi(e_{23})& 0&  0& \varphi(e_{26})\\
          0&  \varphi(e_{32})&               0&       \varphi(e_{34})&   0&      0\\
            0&     0&   \varphi(e_{43})&  0&  \varphi(e_{45})&  0\\
            \varphi(e_{51})&  0&  0&  \varphi(e_{54})&  0&  \varphi(e_{56})\\
            0& \varphi(e_{62})& 0&  0& \varphi(e_{65})&  0
  \end{array}
\right).
\end{displaymath}
By simple calculations, we have

\begin{displaymath}
r(G_{10}^{\varphi})=r\left(
  \begin{array}{cccccccccccccc}
            0&         \varphi(e_{12}) \\
            \varphi(e_{21})&          0 \\
  \end{array}
\right)+r\left(
  \begin{array}{cccccccccccccc}
           0&         \varphi(e_{34}) \\
           \varphi(e_{43})&                    0 \\
  \end{array}
\right)+r\left(
  \begin{array}{cccccccccccccc}
           a+\overline{a}&         \overline{b} \\
           b&                   0 \\
  \end{array}
\right),
\end{displaymath}

where $a=\frac{\varphi(e_{15})\varphi(e_{32})\varphi(e_{54})}{\varphi(e_{12})\varphi(e_{34})}$, $b=\varphi(e_{65})-\frac{\varphi(e_{15})\varphi(e_{62})}{\varphi(e_{12})}$.

So, $r(G_{10}^{\varphi})=4$ if and only if $Re(a)=0$ and $b=0$, i.e., the subgraph induced on vertices $1,~2,~3,~4,~5$ is of Type E and the subgraph induced on vertices $1,~2,~6,~5$ is of Type A.

\begin{displaymath}
A(G_{11}^{\varphi})=\left(
  \begin{array}{cccccccccccccc}
            0&        \varphi(e_{12})&      0& 0&     0&    \varphi(e_{16})&  0\\
     \varphi(e_{21})&               0&        \varphi(e_{23})&     0&  0& 0&   \varphi(e_{27})\\
         0&   \varphi(e_{32})&        0&              \varphi(e_{34})&   0& 0&     0\\
        0&  0&  \varphi(e_{43})&   0&        \varphi(e_{45})&  0&    0\\
          0&  0&  0&  \varphi(e_{54})&  0&    \varphi(e_{56})&  0\\
           \varphi(e_{61})&   0&  0&   0&    \varphi(e_{65})&   0&   \varphi(e_{67})\\
           0&    \varphi(e_{72})&   0&   0&   0&   \varphi(e_{76})&  0
  \end{array}
\right).
\end{displaymath}
By simple calculations, we have
\begin{displaymath}
r(G_{11}^{\varphi})=r\left(
  \begin{array}{cccccccccccccc}
            0&         \varphi(e_{12}) \\
            \varphi(e_{21})&          0 \\
  \end{array}
\right)+r\left(
  \begin{array}{cccccccccccccc}
            0&         \varphi(e_{34}) \\
            \varphi(e_{43})&          0 \\
  \end{array}
\right)+r\left(
  \begin{array}{cccccccccccccc}
               0&         a&                  0\\
    \overline{a}&         0&        \overline{b}\\
               0&         b&                   0
  \end{array}
\right),
\end{displaymath}

where $a=\varphi(e_{56})+\frac{\varphi(e_{16})\varphi(e_{32})\varphi(e_{54})}{\varphi(e_{12})\varphi(e_{34})}$, $b=\varphi(e_{76})-\varphi(e_{16})\frac{\varphi(e_{72})}{\varphi(e_{12})}$.

So, we have $r(G_{11}^{\varphi})=4$ if and only if $a=0,~b=0$.

i.e., $r(G_{11}^{\varphi})=4$ if and only if the subgraphs induced on vertices $1,~2,~3,~4,~5,~6$ and $1,~2,~7,~6$ are of Type A.

If $p\geq0,~l\geq1,~q\geq4$ or $p\geq0,~l\geq2,~q\geq3$ or $p\geq1,~l\geq1,~q\geq4$, by Lemmas \ref{le:2.1}, \ref{le:2.2} and \ref{le:2.3}, then we have $r(\Phi)\geq6$.

This proof is complete.\quad $\square$

\begin{figure}[htbp]
\begin{minipage}[hbt]{1\columnwidth}
\textbf{Table 2}

 {\small {The gain conditions for each gain graph in Theorem \ref{th:4.2} satisfying $r(G^{\varphi})=4$.}}
  \vskip1mm
 {\footnotesize
\begin{tabular}{p{4cm}p{8.8cm}}
\hline
Gain graphs $G^{\varphi}$ & Gain conditions of $G^{\varphi}$\\
\hline
$G^{\varphi}_{12}$ & the subgraph induced on vertices $1,~2,~3$ is of Type E. \\

$G^{\varphi}_{13},~G^{\varphi}_{14},~G^{\varphi}_{19},~G^{\varphi}_{20}$ & any gain. \\

$G^{\varphi}_{15},~G^{\varphi}_{16}$ & the subgraph induced on vertices $1,~2,~3$ is of Type E and the subgraph induced on vertices $1,~2,~4,~3$ is of Type A. \\

$G^{\varphi}_{17},~G^{\varphi}_{18}$ & the subgraph induced on vertices $1,~2,~3,~4$ is of Type A. \\

$G^{\varphi}_{21},~G^{\varphi}_{22}$ & the subgraphs induced on vertices $1,~2,~3,~4$ and $1,~2,~5,~4$ are of Type A. \\

\hline
\end{tabular} }
\end{minipage}
\end{figure}

\begin{figure}[htbp]
  \centering
  \includegraphics[scale=0.6]{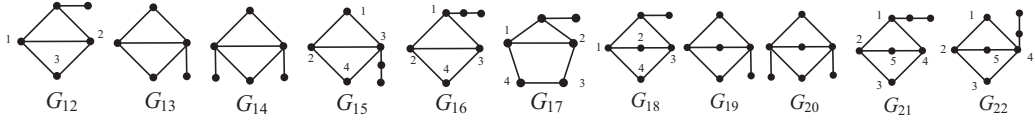}
\caption{The eleven graphs in Theorem 4.2.}
\end{figure}

\noindent\begin{theorem}\label{th:4.2}

Let $\Phi=(G,\varphi)$ be a $\mathbb{T}$-gain bicyclic graph with pendant vertices but no pendant twins. Then $r(\Phi)=4$ if and only if $\Phi$ with one of $G_{i}$'s $(i=12-22)$ (as shown in Fig. 4) as its underlying graph, and the corresponding gain conditions are as shown in Table 2.

\end{theorem}

\noindent\textbf{Proof.} By Theorems \ref{th:3.2}, \ref{th:3.3} and \ref{th:3.4}, we need only characterize all $\mathbb{T}$-gain bicyclic graphs with $\theta(0,1,1)$, $\theta(0,1,2)$, and $\theta(1,1,1)$ as bases.

\textbf{Case 1.} $G$ is a bicyclic graph with $\theta(0,1,1)$ as a base.

\textbf{Subcase 1.1} $G-\theta(0,1,1)$ is a collection of isolated vertices.

   If $\mid G-\theta(0,1,1)\mid=1$, then $G$ is $G_{12}$ or $G_{13}$ (as shown in Fig. 4). For $G_{12}^{\varphi}$, by Lemmas \ref{le:2.2} and \ref{le:2.3}, we have $r(G_{12}^{\varphi})=4$ if and only if the subgraph induced on vertices $1,~2,~3$ in $G_{12}^{\varphi}$ is of Type E.  For $G_{13}^{\varphi}$, by Lemma \ref{le:2.2}, we have $r(G_{13}^{\varphi})=4$, where each edge in $G_{13}^{\varphi}$ has any gain.

    If $\mid G-\theta(0,1,1)\mid=2$, then $G$ is $G_{14}$ (as shown in Fig. 4). By Lemma \ref{le:2.2}, we have $r(G_{14}^{\varphi})=4$, where each edge in $G_{14}^{\varphi}$ has any gain.

   If $\mid G-\theta(0,1,1)\mid\geq3$, by Lemma \ref{le:2.2}, then we have $r(\Phi)\geq6$.

\textbf{Subcase 1.2} If $G-\theta(0,1,1)=P_{2}$, then $G$ is $G_{15}$ or $G_{16}$ (as shown in Fig. 4). By Lemma \ref{le:2.2} and Theorem \ref{th:4.1}, we have $r(\Phi)=4$ if and only if the subgraph induced on vertices $1,~2,~3$ is of Type E and the subgraph induced on vertices $1,~2,~4,~3$ is of Type A in $G_{15}^{\varphi}$ or $G_{16}^{\varphi}$.

\textbf{Subcase 1.3} If $G-\theta(0,1,1)$ contains the union of $P_{2}$ and isolated vertices or contains $P_{3}$ as an induced subgraph, by Lemmas \ref{le:2.2} and \ref{le:2.3}, then we have $r(\Phi)\geq6$.

\textbf{Case 2.} $G$ is a bicyclic graph with $\theta(0,1,2)$ as a base.

\textbf{Subcase 2.1} $G-\theta(0,1,2)$ is a collection of isolated vertices.

   If $\mid G-\theta(0,1,2)\mid=1$, then $G$ is $G_{17}$ (as shown in Fig. 4). By Lemmas \ref{le:2.2} and \ref{le:2.3}, we have $r(G_{17}^{\varphi})=4$ if and only if the subgraph induced on vertices $1,~2,~3,~4$ in $G_{17}^{\varphi}$ is of Type A.

  If $\mid G-\theta(0,1,2)\mid\geq2$, by Lemma \ref{le:2.2}, then we have $r(\Phi)\geq6$.

\textbf{Subcase 2.2} $G-\theta(0,1,2)$ contains $P_{2}$ as an induced subgraph.

In this case, by Lemma \ref{le:2.2} and Theorem \ref{th:4.1}, we have $r(\Phi)\geq6$.

\textbf{Case 3.} $G$ is a bicyclic graph with $\theta(1,1,1)$ as a base.

\textbf{Subcase 3.1} $G-\theta(1,1,1)$ is a collection of isolated vertices.

 If $\mid G-\theta(1,1,1)\mid=1$, then $G$ is $G_{18}$ or $G_{19}$ (as shown in Fig. 4). For $G_{18}^{\varphi}$,  by Lemmas \ref{le:2.2} and \ref{le:2.3}, we have $r(G_{18}^{\varphi})=4$ if and only if the subgraph induced on vertices $1,~2,~3,~4$ in $G_{18}^{\varphi}$ is of Type A. For $G_{19}^{\varphi}$, by Lemma \ref{le:2.2},  we have $r(G_{19}^{\varphi})=4$, where each edge in $G_{19}^{\varphi}$ has any gain.

 If $\mid G-\theta(1,1,1)\mid=2$, then $G$ is $G_{20}$ (as shown in Fig. 4). By Lemma \ref{le:2.2}, we have $r(G_{20}^{\varphi})=4$, where each edge in $G_{20}^{\varphi}$ has any gain.

 If $\mid G-\theta(1,1,1)\mid\geq3$, by Lemma \ref{le:2.2}, then we have $r(\Phi)\geq6$.

\textbf{Subcase 3.2} If $G-\theta(1,1,1)=P_{2}$, then $G$ is $G_{21}$ or $G_{22}$ (as shown in Fig. 4). By Lemma \ref{le:2.2} and Theorem \ref{th:4.1}, we have $r(\Phi)=4$ if and only if the subgraphs induced on vertices $1,~2,~3,~4$ and $1,~2,~5,~4$ are of Type A in $G_{21}^{\varphi}$ or $G_{22}^{\varphi}$.

\textbf{Subcase 3.3} If $G-\theta(1,1,1)$ contains the union of $P_{2}$ and isolated vertices or contains $P_{3}$ as an induced subgraph, by Lemmas \ref{le:2.2} and \ref{le:2.3}, then we have $r(\Phi)\geq6$.

This proof is complete.\quad $\square$

\end{document}